\documentclass[11pt,a4paper]{article}
\usepackage[utf8]{inputenc}
\usepackage{mathtools}

\usepackage[numbers]{natbib}

\usepackage{amsthm}
\usepackage{array}
\usepackage{multirow}

\usepackage{epigraph} 
 
\usepackage{setspace} 
 \usepackage{hyperref}
 \usepackage{amssymb}
\usepackage{amsfonts}

\usepackage{caption} 
\usepackage{needspace} 
 \usepackage{centernot}

\usepackage{bussproofs}

\usepackage{proof}

\usepackage[tableau]{prooftrees}

\usepackage{tikz}
\usetikzlibrary{matrix}

\usepackage[switch]{lineno}

\title{\bf Analytic proofs for logics of evidence and truth  
 \thanks{The research of the first author is supported by the grant 307376/2018-4, 
Conselho Nacional de Desenvolvimento Cient\'ifico e Tecnol\'ogico (CNPq, Brazil). 
 The research of the third author is supported by the grants 311911/2018-8,
Conselho Nacional de Desenvolvimento Cient\'ifico e Tecnol\'ogico (CNPq, Brazil), and APQ-02093-21,
Funda\cao\ de Amparo  \`a Pesquisa do Estado de Minas Gerais (FAPEMIG, Brazil).  
 }}

\author{Walter Carnielli$^1$, Lorenzzo Frade$^2$, Abilio Rodrigues$^2$  \\
[4mm]
$^1$University of Campinas, Campinas, Brazil \\
walterac@unicamp.br \\
$^2$Federal University of Minas Gerais, Belo Horizonte, Brazil \\
lorenzzo\_frade@hotmail.com  \\ 
abilio.rodrigues@gmail.com 
}

\usepackage[shortlabels]{enumitem}

\label{xEXPRESSIONS}

\label{xLOGICS}{}

\newcommand{\lfis}{\textit{LFI}s}
\newcommand{\lfi}{\textit{LFI}}

\newcommand{\letf}{$LET_{F}$}

\newcommand{\fde}{$FDE$}

\newcommand{\lets}{\textit{LET}s}

\label{xCOLORS} %

\label{xSYMBOLS}

\newcommand{\cons}{\ensuremath{{\circ}}} 
\newcommand{\incon}{\ensuremath{{\bullet}}}

\newcommand{\wneg}{\ensuremath{\lnot}}
\newcommand{\sneg}{\ensuremath{{\sim}}}

\label{xGENERAL COMMANDS}%

\newcommand{\mh}{\noindent}

\newcommand{\m}{\vspace{1mm}}
\newcommand{\mm}{\vspace{2mm}}
\newcommand{\mmm}{\vspace{3mm}}
\newcommand{\mmmm}{\vspace{4mm}}
\newcommand{\setl}{\setlength\itemsep{-0.2em}}

\newcommand{\bqu}{\begin{quote} \small}
\newcommand{\equ}{\end{quote} \normalsize}
\newcommand{\enr}{\begin{enumerate}[label={(\arabic*)}, resume]}
\newcommand{\eenr}{\end{enumerate}}

\usepackage{graphicx}
\makeatletter
\providecommand*{\ndashv}{%
  \mathrel{%
    \mathpalette\@ndashv\nvdash
  }%
}
\newcommand*{\@ndashv}[2]{%
  \reflectbox{$\m@th#1#2$}%
}
\makeatother

\label{xPORTUGUES}

\newcommand{\cao}{\c{c}\~{a}o}

\theoremstyle{definition}
\newtheorem{Th}{Theorem}
\newtheorem{Ex}[Th]{Example}

\newtheorem{Le}[Th]{Lemma}
\newtheorem{Df}[Th]{Definition}

\date{}

\begin{document}

\maketitle
    

\abstract{\mh 
This paper presents  a  sound, complete, and decidable analytic
tableau  system for  the logic of evidence and truth \letf, introduced in  
 Rodrigues, Bueno-Soler \& Carnielli (Synthese, DOI: 10.1007/s11229-020-02571-w, 2020).   
\letf\ is an extension of the logic of
first-degree entailment (\fde), also known as Belnap-Dunn logic. \fde\ is
a  widely studied four-valued  paraconsistent logic, 
with applications in computer
science  and in the algebra of processes. 
\letf\  extends \fde\  in a very 
natural way,   by adding  a classicality operator \cons, which recovers 
classical logic for propositions in its scope, 
and  a non-classicality operator  \incon,  dual of \cons.   
}


\mmmm

\begin{flushright}
\textit{Dedicated to Francisco Mir\'o Quesada, \\ the godfather of paraconsistency.}
\end{flushright}

\section{Introduction}

 The aim of this paper is to present an analytic tableau system for the  logic \letf, which  is a member of a family of logics called \textit{logics of evidence and truth} (\textit{LET}s) \cite{letj, letf}.  
The motivation  for \lets\ is the idea that in real-life reasoning we deal with positive and negative
evidence, and such evidence can be  conclusive or non-conclusive.  
  \lets\ thus   combine 
 two different notions of logical consequence  in the same formal system: one   
preserves truth  (classical consequence), the other   preserves evidence -- hence, the name `logics of evidence and truth'.  
Evidence is thought of as a notion weaker
 than truth, in the sense that there may be evidence for a proposition $A$ even in the case  $A$ is not true.  
Positive and negative evidence, 
respectively evidence for truth and for falsity, are independent
and non-complementary, 
and negative evidence for  $A$ is identified with positive evidence for
$\neg A$. 
\lets\ are paraconsistent, since it may be that there is conflicting non-conclusive 
evidence for a proposition $A$, and paracomplete, 
since  it may happen that there is no evidence at all for  $A$.\footnote{A more detailed account of the notion of evidence that 
underlies the intended interpretation of \lets\ can be found in \cite{answer.barrio}. }

The logic \letf,  introduced in \cite{letf}, 
is an extension of the logic of first-degree entailment (\fde), also called Belnap-Dunn   logic 
\cite{belnap1977.how, dunn76,essays.belnap.dunn}.  
\letf\  is equipped with a classicality operator \cons\ and a non-classicality operator 
\incon\  which is the dual  of \cons. The deductive behavior of  \cons\ and \incon\ is given by the following inferences:
\enr \setl 
\item  $ \cons A, \incon A   \vdash  B$,  
\item   $\vdash \cons A \lor \incon A $, 
\item  $\cons A, A, \neg A \vdash B$  (although  $A, \neg A \nvdash  B$), \label{gent.exp}
\item  $\cons A \vdash A \lor \neg A$ (although $\nvdash A \lor \neg A $), \label{gent.pem}
\item   $\vdash A \lor \neg A \lor \incon A $,  \label{dual.gent.exp}
\item   $A , \neg A \vdash \incon  A $. \label{dual.gent.pen}
\eenr 

\mh According to the intended interpretation in terms of evidence, 
items (1) and (2) together mean that either there is or there is not conclusive evidence for $A$. Items \ref{gent.exp} and 
\ref{gent.pem}, as in any \textit{LET}, recover classical negation for propositions in the scope of \cons, and 
\ref{dual.gent.exp} and \ref{dual.gent.pen} are dual, respectively, of \ref{gent.exp} and 
\ref{gent.pem}. 
 Due to items \ref{gent.exp} and \ref{gent.pem},  \letf\ is a  logic  of formal inconsistency and undeterminedness 
 \cite[cf.][]{recovery, marcos2005}.  
 
 \letf, as \lets\ in general, can also be interpreted in terms of information, which may be   unreliable or  reliable, the latter being 
subjected to classical logic \cite[see][Sect.~2.2.1]{letf}. In this case, $\cons A$ means that the information about $A$, positive or negative, is reliable. In \cite{kpk.let}, Kripke models for \letf\ have been proposed. 
These models intend to represent a database that
receives information as time passes, and such information can be positive, negative, unreliable,
or reliable. This idea fits  the interpretation of  Belnap-Dunn   logic  as an information-based
logic, but adds to the four scenarios expressed by it two new scenarios: reliable  
information  (i) for the truth and (ii) for the falsity of a given proposition. 
\fde\  was probably the first  logic  to be  applied in 
computer science and, more recently, in description  logics  and in the
algebra of processes   \cite[see e.g.][]{bergstra,ginsberg,patel,straccia}. 
An extension of \fde\ such as  \letf\  has an ample range of possible 
applications, especially in Bayesian decision   procedures under
uncertainty  \cite[cf.][]{future}. 

This paper extends the investigation on \letf\ carried out in \cite{kpk.let,letf}. It has been remarked in \cite[Sect.~3.2]{letf} that \letf\ is decidable, but a detailed algorithm has not been presented. 
Our main aim here is to present a sound, complete, and decidable  tableaux system for \letf.\footnote{Analytic tableaux for some logics of formal inconsistency have already been proposed in \cite{tableau.lfi} and \cite[Sect.~3.5]{carcoma2007}.}  
The tableaux to be proposed   are analytic because the formulas yielded by the application of a rule to a formula \textit{F} 
are always less complex than  \textit{F}, and counterexamples can be obtained from open branches of terminated tableaux. 

The remainder of this paper is structured as follows. Section \ref{sec.val.sem} 
presents the valuation semantics of \letf, and   Section  \ref{sec.tab}, the corresponding  tableau  system. 
In Section~\ref{completeness} we prove the soundness and completeness of the \letf-tableau system with respect 
to  the semantics of Section~\ref{sec.val.sem}, and   also  that the analytic 
 tableau system provides a decision procedure for \letf. 
Finally, in Section  \ref{examples},  some examples of \letf-tableaux are given and commented.

\section{Valuation semantics for \letf} \label{sec.val.sem}

The language $\mathcal{L}$ of \letf\ is composed of denumerably many sentential letters $p_{1}, p_{2}, \dots$, the~unary connectives $\cons$, \incon, and $\neg$, the~binary connectives $\land$ and $\lor$, and~parentheses. The~set of formulas of $\mathcal{L}$, 
which is also denoted by $\mathcal{L}$, is inductively defined in the usual way. 
 Roman capitals $A, B, C,\dots$ will be used as metavariables for the formulas of $\mathcal{L}$,  and Greek capitals 
$\Gamma, \Delta, \Sigma, \dots$  as metavariables for sets of formulas $\mathcal{L}$. 

\m 

In \cite{letf} a natural deduction system, together with the following sound and complete semantics,  was presented for \letf. 

\begin{Df}  {[Valuation semantics for \letf]}  \label{val.sem.letf}

\mh
A  valuation semantics for \letf\  is a collection of \letf-valuations defined as follows. 
A function $v:  \mathcal{L} \to \{0,1\}$ is a \letf-valuation  if it  satisfies the following clauses:
\begin{itemize} \setl 
\item[(v1)] $v(A \land  B )=1$ iff $v(A)=1$ and $v( B )=1$,
\item[(v2)] $v(A \lor  B )=1$ iff $v(A)=1$ or $v( B )=1$,
\item[(v3)] $v(\neg (A \land  B ))$ = 1 iff $v(\neg A) = 1$ or  $v(\neg  B ) = 1$,
\item[(v4)] $v(\neg (A \lor  B )) = 1$ iff $v(\neg A) = 1$ and  $v(\neg  B ) = 1$,
\item[(v5)] $v(A) = 1$ iff $v(\neg\neg A) = 1$, 
\item[(v6)] If $v(\cons A) = 1$, then  $v(A) = 1$ if and only if $v(\neg A) = 0$,
\item[(v7)] $ v(\incon A) = 1$ iff $ v(\cons A) = 0$.  
\end{itemize}
\end{Df}

\begin{Df}
We say that a formula $ A$ is a semantical consequence of $\Gamma$, $\Gamma \vDash  A$, 
if and only if for 
every valuation $v$, if $v( B )=1$ for all $ B \in \Gamma$, then $v(A) = 1$.  
\end{Df}

\mh The semantics above is non-deterministic in the sense that the semantic value of  complex formulas  
is not always functionally determined by its parts, as the following non-deterministic 
matrix (also called  quasi-matrix)  shows.\footnote{On non-deterministic valuation semantics, see Loparic \cite[][]{loparic1986,loparic2010}.}

\begin{Ex} \ \label{ex.valuations.letf}

\m 
\mh \quad \quad (i) $\cons p \vDash p \lor \neg p$,  

\mh \quad \quad  (ii) $p \lor \neg p \nvDash \cons p$, 

\mh\quad\quad (iii) $\nvDash \cons p \lor \neg \cons p$.

\mmm 
\mh 
\begin{scriptsize}
\begin{tabular}{|c|c|c|c|c|c|c|c|c|c|c|c|c|c|}
\hline 
$_1$ &  $p$ & \multicolumn{6}{c|}{0} & \multicolumn{6}{c|}{1}\tabularnewline
\hline 
$_2$ &  $\neg p$ & \multicolumn{2}{c|}{0} & \multicolumn{4}{c|}{1} & \multicolumn{4}{c|}{0} & \multicolumn{2}{c|}{1}\tabularnewline
\hline 
$_3$ &  $p \lor\neg p$ & \multicolumn{2}{c|}{0} & \multicolumn{4}{c|}{1} & \multicolumn{4}{c|}{1} & \multicolumn{2}{c|}{1}\tabularnewline
\hline 
$_4$ & $\cons p$ & \multicolumn{2}{c|}{0} & \multicolumn{2}{c|}{0} & \multicolumn{2}{c|}{1} & \multicolumn{2}{c|}{0} & \multicolumn{2}{c|}{1} & \multicolumn{2}{c|}{0}\tabularnewline
\hline 
$_5$ & $\neg\cons p$ & 0 & 1 & 0 & 1 & 0 & 1 & 0 & 1 & 0 & 1 & 0 & 1\tabularnewline
\hline 
 &  & $v_{1}$ & $v_{2}$ & $v_{3}$ & $v_{4}$ & $v_{5}$ & $v_{6}$ & $v_{7}$ & $v_{8}$ &$v_{9}$  &$v_{10}$  &$v_{11}$  & $v_{12}$\tabularnewline
\hline 
\end{tabular}
\end{scriptsize}

 \mmm \mh There is no valuation $v$  such that $v(\cons p) = 1$ and $v(p \lor\neg p) = 0$, so (i) holds. Valuation 
  $v_3$  provides a counterexample  to both (ii) and (iii). Note that lines 2,   4, and 5 bifurcate: line 2 because $p$ 
  and $\neg p$ are completely independent of each order, line 4 because there are no sufficient 
  conditions  
  for $v(\cons A) = 1$,  and line 5 for the same reason as line 2. 
Indeed, a   feature of \letf\ (as well as of some \lfis) is that when $\cons A$ ($\incon A$) occurs in the scope of $\neg$, 
unless $\cons\cons A$ ($\cons\incon A$) holds, the value of 
$\neg\cons A$ ($\neg\incon A$) is not functionally determined by the value of $\cons A$ ($\incon A$). 
The negation in these formulas is still a weak negation. 
  
\end{Ex}

\section{An analytic tableau system for \letf} \label{sec.tab}

Given the semantics above, 
we shall prove in Section \ref{completeness}   that the following tableau rules
constitute a sound, complete, and decidable proof system for \letf. 
We will consider informally here the usual notions related to tableaux: trees,
branches, nodes, etc. The
labels 0 and  1  refer to metamathematical  markers, intuitively
related to the semantic values \textit{0} and \textit{1}.


\begin{Df} {[Tableau rules for \letf]}

\m 

 {
	\needspace{2\baselineskip}
	\label{tab:regras-tablos}

	\begin{tableau}
				{		to prove={\textrm{\textbf{Rule 1}}}, 
		for tree={s sep=10mm},line numbering=false 
	}
	[1(A \land B)
	[1(A)
	[1(B)
	]
	] 
	]
	\end{tableau}
	\qquad  \qquad  \qquad  \qquad  \qquad  \qquad \, \, \,
	\begin{tableau}
		{
	to prove={\textrm{\textbf{Rule 2}}}, 
	for tree={s sep=10mm},line numbering=false 
}
[0(A \land B)
[0(A)]
[0(B)] 
]
	\end{tableau}
	
	\bigskip

	\begin{tableau}
	{
		to prove={\textrm{\textbf{Rule 3}}}, 
		for tree={s sep=10mm},line numbering=false 
	}
	[1(\neg(A \land B))
	[1(\neg A)]
	[1(\neg B)]
	]
\end{tableau}
\qquad  \qquad  \qquad  \qquad   \qquad 
\begin{tableau}
	{
		to prove={\textrm{\textbf{Rule 4}}}, 
		for tree={s sep=10mm},line numbering=false 
	}
	[0(\neg(A \land B))
	[0(\neg A)
	[0(\neg B) 
	]
	]
	]
\end{tableau}

\bigskip

	\begin{tableau}
	{
		to prove={\textrm{\textbf{Rule 5}}}, 
		for tree={s sep=10mm},line numbering=false 
	}
	[1(A \lor B)
	[1(A)]
	[1(B)]
	]
\end{tableau}
\qquad  \qquad  \qquad  \qquad  \qquad  \qquad 
\begin{tableau}
	{
		to prove={\textrm{\textbf{Rule 6}}}, 
		for tree={s sep=10mm},line numbering=false 
	}
	[0(A \lor B)
	[0(A)
	[0(B) 
	]
	]
	]
\end{tableau}

\bigskip

\begin{tableau}
	{
		to prove={\textrm{\textbf{Rule 7}}}, 
		for tree={s sep=10mm},line numbering=false 
	}
	[1(\neg(A \lor B))
	[1(\neg A)
	[1(\neg B)
	]
	]
	]
\end{tableau}
\qquad  \qquad  \qquad  \qquad  \qquad  \qquad
\begin{tableau}
	{
		to prove={\textrm{\textbf{Rule 8}}}, 
		for tree={s sep=10mm},line numbering=false 
	}
	[0(\neg(A \lor B))
	[0(\neg A)]
	[0(\neg B)] 
	]
\end{tableau}

\bigskip

\begin{tableau}
	{
		to prove={\textrm{\textbf{Rule 9}}}, 
		for tree={s sep=10mm},line numbering=false 
	}
	[1(\neg\neg A)
	[1(A)
	]
	]
\end{tableau}
\qquad  \qquad  \qquad  \qquad  \qquad  \qquad  \qquad  \qquad 
\begin{tableau}
	{
		to prove={\textrm{\textbf{Rule 10}}}, 
		for tree={s sep=10mm},line numbering=false 
	}
	[0(\neg\neg A)
	[0(A)
	]
	]
\end{tableau}

\bigskip

\begin{tableau}
	{
		to prove={\textrm{\textbf{Rule 11}}}, 
		for tree={s sep=10mm},line numbering=false 
	}
	[1(\cons A)
	[1(A)
	[0(\neg A)]]
	[0(A)
	[1(\neg  A)]]
	]
\end{tableau}

\bigskip

\begin{tableau}
	{
		to prove={\textrm{\textbf{Rule 12}}}, 
		for tree={s sep=10mm},line numbering=false 
	}
	[1(\incon A)
	[0(\cons A)
	]
	]
\end{tableau}
\qquad  \qquad  \qquad  \qquad  \qquad  \qquad  \qquad  \qquad 
\begin{tableau}
	{
		to prove={\textrm{\textbf{Rule 13}}}, 
		for tree={s sep=10mm},line numbering=false 
	}
	[0(\incon A)
	[1(\cons A)
	]
	]
\end{tableau}

}

\end{Df}
 
\mm

 There is no need for a rule for $0(\cons A)$. Such a rule, call it \textbf{R}, would conclude $1(\incon A)$ from $0(\cons A)$. 
Besides yielding a loop with \textbf{Rule 12}, it can be shown that  
\textbf{R} is not necessary at all.   
Suppose the application of \textbf{R} to $0(\cons A)$ yielded a closed  branch $b$  
such that both $1(\incon A)$ and $0(\incon A)$ occur in $b$. But in this case, it would be enough to apply \textbf{Rule 13} 
to $0(\incon A)$, obtaining  a branch $b'$ containing $1(\cons A)$ and $0(\cons A)$, 
and $b'$ 
would be a closed branch.  

Concerning  \textbf{Rule 11}, recall that the symbol \cons\ in \letf\ expresses  \textit{classicality}, i.e., 
a formula $\cons A$ implies that $A$ 
behave  classically. 
The classical behavior of $A$ is recovered by recovering classical negation for $A$:     
either $A$, or $\neg A$ holds, and not both. This is precisely what \textbf{Rule 11} does. 
The semantic clause for $\cons A$ has only a necessary condition for 
$v(\cons A) = 1$, and the absence of a rule for $0(\cons A)$ mimics this fact.

Moreover, note that 
there are no tableau rules for $\neg\cons A$  and $\neg\incon A$.
As remarked  in Example \ref{ex.valuations.letf} above,  
the  value of 
$\neg\cons A$ ($\neg\incon A$) is not   functionally determined by the value of $\cons A$ ($\incon A$), except when 
$\cons\cons A$ ($\cons\incon A$) holds. 
This is illustrated by the fact that in \letf\ 
neither $\vdash \cons A \lor \neg\cons A$ , nor  $\cons A, \neg\cons A \vdash B$ hold (we will see the respective 
 tableaux in Section \ref{examples}  below).

\m

\begin{Df} {[\letf\ tableaux]} \label{varios} 
\begin{enumerate} \setl 
\item We define a \textit{tableau} for a set $\Delta$ of signed formulas as a
tree whose  first node  contains all the signed formulas in $\Delta$, and whose
subsequent   nodes are obtained  by applications of the tableau rules given in Definition \ref{tab:regras-tablos} above.

 \item A tableau branch is \emph{closed} if it contains a  pair  of  signed formulas $1(F)$ and $0(F)$. 
    If a  branch is not closed, we say it is  open. 
   A tableau is
{closed} if all its branches are closed.

\item
 When no rule can be applied to any open branch, the 
 tableau is \textit{terminated}. Every closed  tableau is also a terminated tableau.

\item
A  closed tableau is \emph{pruned} if all formulas in $\Delta$ that have not been
used in its branches are  deleted.

\item
A formula   $A$ has a \textit{proof}  from premisses $\Gamma$, denoted
$\Gamma\vdash A$,   if there  exists a closed  tableau  for the set
$\{ 1(B): B \in \Gamma \} \cup \{ 0(A) \}$.

$\Gamma$ can be empty, and in this case a  
proof of $ \vdash A$ 
reduces to a tableau for   the singleton $\{ 0(A) \}$.
 
\end{enumerate}

\end{Df}

\m 

\begin{Th} {[Compacteness of tableau proofs]} \ \label{compactness}

\m\mh    $\Gamma\vdash A$ iff   $\Gamma_0\vdash A$, for     $\Gamma_0
\subset \Gamma$,    $\Gamma_0$ finite.
\begin{proof}
First, note  that a  pruned closed tableau  is  finite, since
any tableau is a finite collection of  branches and  every closed
branch is  finite. So, after a  pruning procedure, only a finite  subset
$\Gamma_0 \subset \Gamma$  has  been used  (that is, the pruning
procedure deletes  any  possibly  infinite collection of sentences of
$\Gamma$  not used in the process of tableau closure).  
\end{proof}
\end{Th}

\m

An analytic tableau is a procedure of reductio ad absurdum 
whose aim is to obtain a closed tableau, and thus proving $\Gamma\vdash A$. 
A closed  tableau is a halting condition of a  decision procedure, and in the case of \letf\ this halting condition 
is always obtained, since every \letf-tableau terminates (cf.~Theorem \ref{terminates} below). 
If the tableau is open, an open branch gives a valuation 
that is a counterexample for $\Gamma\vdash A$.  We will see that  \letf-tableaux constitute an elegant decision procedure for \letf.

\mm

\section{Soundness and completeness} \label{completeness}

In this section we prove soundness, completeness, and decidability of \letf-tableaux. 
We start by defining the complexity $\mathcal{C}$ of a formula in the language of \letf.

  \begin{Df} {[Complexity]} \label{def-complexity}
\mh 
	The \textit{complexity} of a formula $F$ of the language of \letf\ is given by the map $ \mathcal{C} : F  \rightarrow   \mathbb{N} $, such that:
	\begin{enumerate} \setl 
		\item $ \mathcal{C}(p) = 0$;
		\item  $\mathcal{C}(\neg A) =  \mathcal{C}(A)  + 1$;
		\item   
		$\mathcal{C}(A \ast B) =   \mathcal{C}(A) +  \mathcal{C}(B)  + 1$, 
	 for $\ast \in \{ \land, \lor \}$; 

		\item $ \mathcal{C}(\cons A) =  \mathcal{C}(A)  + 2$; \label{cons.comp}
		\item  $ \mathcal{C}(\incon A) =  \mathcal{C}(A ) + 3$. \label{incon.comp} 

	\end{enumerate}	 

	\end{Df}

\subsection{Soundness}

\begin{Th} {[Termination]} \label{terminates}

		\mh Every \letf-tableau terminates: after a 
	finite number of steps no more  rules can be applied.
	\begin{proof} 
		The result follows from the fact that each rule results in formulas with  less complexity or formulas to which no rule 
		can be applied.  
	\end{proof}
\end{Th}

\m

\begin{Df} {[\letf-satisfiable branch]} \label{letf.sat.branch}
	
	\mh A branch $b$ is  \letf-\textit{satisfiable} if  
		there is a \letf-valuation $v$ such that for every formula that occurs in $b$: 

\mh \quad \quad 1. If $1(F)$ is in $b$, then $v(F) = 1$,

\mh \quad \quad 2. If $0(F)$ is in $b$, then $v(F) = 0$.

\mh In this case, we say that the valuation $v$ satisfies the branch $b$. 

\end{Df}

\mh Clearly, if a branch is closed, it cannot be satisfiable, because there is no 
valuation such that $v(F)=1$ and $v(F)=0$. 
It is worth noting  that a \letf-satisfiable branch, as expected,  differs from a satisfiable branch in classical logic. 
A  \letf-branch $b$ may be satisfiable even if $1(F)$ and $1(\neg F)$, as well as $0(F)$ and $0(\neg F)$, 
are both in $b$.  The condition for a branch $b$ of a \letf-tableau to be closed 
(cf.~Definition \ref{varios} item 2) is that for some $F$, $1(F)$ and $0(F)$ are in $b$.
Moreover, although  there is no \letf-valuation such that $v(F) = 1$, $v(\neg F) = 1$, and $v(\cons F)  = 1$, 
and no open terminated branch $b$ can be such that $1(F)$, $1( \neg F)$, and $1( \circ F)$ are all in $b$, we do not need 
to add the latter as a condition for a closed branch. 
If $1(\circ F)$ is in $b$, \textbf{Rule 11} must be applied, and $b$ will bifurcate into two branches: $b'$, with 
$1(F)$ and $0(\neg F)$, and $b''$, with $0(F)$ and $1(\neg F)$. If $1(F)$ and $1( \neg F)$ are both in $b$, then both branches $b'$ and  $b''$ are closed.

\m

\begin{Le} {[Satisfiable branches]}  \label{lema-ramos-correcao}
	
	\mh (i)  If a non-branching rule is
	applied to a \letf-satisfiable branch, the result is another \letf-satisfiable
	branch. 
	
	\mh (ii) If a branching rule is applied to a \letf-satisfiable branch, at
	least one of the resulting branches is also a \letf-satisfiable branch.
	
	\begin{proof} \

\mh 		To prove (i), we have to show that every non-branching rule applied to a  satisfiable branch results in another  satisfiable branch. 

		\mh For \textbf{Rule 1}, suppose $b$ is a satisfiable branch containing  $1(A  \land  B)$.  
		Therefore, by Definition \ref{letf.sat.branch}, there is a  valuation $v$ such that $v$  satisfies $b$ and $v(A \land B)=1$. 
		Now, applying \textbf{Rule 1} to $b$ yields a branch $b'$ such that $1(A)$ and $1(B)$, as well as $1(A  \land  B)$, are  in $b'$. 
By Definition \ref{val.sem.letf}, since $v(A \land B)=1$, we have that $v(A) =1$ and $ v(B) = 1$. 
		Therefore, $b'$ is 	a satisfiable branch. 
		Analogous reasoning applies to \textbf{Rule 4}, \textbf{Rule 6}, \textbf{Rule 7}, \textbf{Rule 9}, and \textbf{Rule 10}. 
\mh For  \textbf{Rule 12}, 
suppose $b$ is a satisfiable branch containing $1(\incon A)$.  
By Definition \ref{letf.sat.branch}, there is a valuation $v$ such that $v(\incon A)=1$,    and by Definition \ref{val.sem.letf}, 
$v(\cons A) = 0$.  
An application of \textbf{Rule 12} yields a branch $b'$ such that $0(\cons A)$ is in $b'$, and $v$ satisfies $b'$. 
Analogous reasoning applies to \textbf{Rule 13}.
		
		\m 
	\mh	To prove (ii), we must show that every branching rule applied to a satisfiable branch results in at least one satisfiable branch.

	\mh		For \textbf{Rule 2}. Suppose $b$ is a satisfiable branch containing   $0(A  \land  B)$.  
Thus there is a valuation  $v$ such that  $v$ satisfies $b$ and either $v(A)=0$   or $v(B)=0$. 
An application of \textbf{Rule 2} yields two branches: $b'$ and $b''$, such that  $0(A)$ occurs in $b'$ and $0(B)$ in $b''$. 
Therefore, either $v$ satisfies $b'$ or $v$ satisfies $b''$. 
Analogous  reasoning applies to \textbf{Rule 3}, \textbf{Rule 5}, and \textbf{Rule 8}. 
	For \textbf{Rule 11}, suppose $b$ is a satisfiable branch containing $1(\cons A)$. 
Thus there is a valuation 
		$v$, such that $v$ satisfies  $b$ and $v(\cons A) = 1$. Then by Definition \ref{val.sem.letf}  either $v(A) = 1$ and $v(\neg A) = 0$, or  		$v(A) = 0$ and $v(\neg A) = 1$. 
Now, if we apply \textbf{Rule 11} to $b$, we  get two branches: $b'$  with $1(A)$ and $0(\neg A)$,  
and $b''$ with $0(A)$ and $ 1(\neg A) $. Therefore, either $v$ satisfies $b'$, or $v$ satisfies $b''$. 
 \end{proof}
\end{Le}

\m

\begin{Th} {[Soundness]} 
 If $\Gamma \vdash A$, then $\Gamma \vDash A$.
	
	\begin{proof}
		Assume $\Gamma \vdash A$ and suppose for reductio  that $\Gamma \nvDash A$. 
		Let $b$  be the first node  of the tableau for $\Gamma \vdash A$, with every formula $F$ from $\Gamma$ labeled with 1 
		($1(F_1)$, $1(F_2)$, etc.) and $0(A)$. 
	It follows from $\Gamma \nvDash A$ that there is a \letf-valuation $v$ such that for every formula $F$ from 
	$\Gamma$, $v(F) = 1$, and~$v(A) = 0$. Therefore, by Definition \ref{letf.sat.branch}, $b$ is 
 	satisfiable. 
	But from Lemma \ref{lema-ramos-correcao} we have that if any rule is applied to $b$, at least one of the resulting branches will be satisfiable. Hence, after a finite number of rule applications, the tableau for $\Gamma \vdash A$ terminates with (at least) one satisfiable branch $b'$. 
	There is thus a valuation $v'$ such that for every formula $F$ that occurs in $b'$, if $1(F)$ is in $b'$, $v(F)=1$, 
	and if $0(F)$ is in $b'$, $v(F)=0$,  
	and 	$b'$ cannot contain any formula $F$ such that $1(F)$  and $0(F)$  are both in $b'$ --  otherwise, 
	$v'$  would not be a valuation. 
But then, $b'$ is  an open branch, and $\Gamma \nvdash A$, which contradicts the initial assumption.  
		Therefore, $\Gamma \vDash A$. 		
	\end{proof}	
\end{Th}

\m

\subsection{Completeness}

In order to prove completeness of \letf-tableaux with respect to the valuation semantics of \letf,  
we  show the contrapositive: if $\Gamma$ $\nvdash A$, then $\Gamma$ $\nvDash A$. 
 $\Gamma$ $\nvdash A$ just in case 
there is an open branch in the tableau. Let $b$ be this open branch. 
We have to show that there is a valuation $v$ induced by $b$  such that for every formula $F$ from $\Gamma$, 
$v(F) = 1$, and $v(A) = 0$. 
We therefore begin by defining a valuation induced by an open branch. 

\begin{Df} {[Semi-valuation induced by an open branch]} \label{semi.val}

	\mh Let a \textit{literal} be  a propositional letter or the negation of a propositional letter. 
	Let $b$ be an open branch of a terminated tableau. The semi-valuation $s$ induced by $b$ is   such that:
	\begin{enumerate}
\setl 
		\item For every literal $l$  such that $1(\textit{l})$ is in $b$, $s(l) = 1$; 
		\item For every literal $l$  such that $1(\textit{l})$ is not in $b$, $s(l) = 0$;  
		\item If $ 1(\cons A) $ is in $ b $, then 
		$ s(\cons A) = 1 $ and $ s(\incon A) = 0 $;  
		\item If $1(\cons A)$ is not in $b$, then $s(\cons A) = 0$ and $s(\incon A) = 1$; 
		\item  $s(\neg \cons  A) = 1$ if, and only if, $ 1(\neg\cons A) $ is in $ b $;  
		\item  $s(\neg\incon A) = 1$ if, and only if, $1(\neg\incon A)$ is in $ b $.
		
\end{enumerate} 	
\end{Df}

\mh  As remarked  in   Example~\ref{ex.valuations.letf}, the semantic values of 
$\neg\cons A$ and $\neg\incon A$ are not  functionally determined by the values of $\cons A $ and $\incon A$. 
It is for this reason that the 
items 5 and 6 in Definition~\ref{induced.val} above have to be explicitly given.  
 
\mm

\begin{Le} {[Valuation induced by an open branch]} \label{induced.val} \ 

\mh 	Let $b$ be an open branch of a terminated tableau. 
Then, there exists a valuation $v$ induced by $b$ such that  for every formula $F$:

 \mh 	\quad	(i) If $1(F)$ is in $b$, $v(F) = 1$, 
		
	\mh	\quad (ii) If $0(F)$ is in $b$, $v(F) = 0$. 
 		\begin{proof}
		\mh

	\mh The proof is by  induction on the complexity of $F$ (Definition~\ref{def-complexity}).  

\m 

\mh 	(1)	If $F$ is a literal,  $\cons A$, $\incon A$, $\neg\cons A$, or $\neg \incon A$, define $v(F) = s(F)$. 
\label{def.val.semi}
		
	\mm 
		
\mh (2)  $F = A \land B$. 

\m 

\mh (2.1) If $1(A \land  B)$ is in $b$,  since the  tableau is terminated, 
 \textbf{Rule 1} has been applied, therefore $1(A)$ is in $b$ and $1(B)$ is in $b$. 
 By inductive hypothesis, $v(A) = 1$ and $v(B) = 1$. We then define  $v(A  \land   B) = 1$.
 
 \m 
 
\mh (2.2)   If $0(A  \land   B)$ is in $b$, \textbf{Rule 2} was applied and the  tableau bifurcated into two branches: 
 $b'$ and $b''$. 
In $b'$, we have $0(A)$, and in $b''$, $0(B)$. 
 Since $b$  is an open branch of a terminated tableau, we have two (non-excluding) options: 
 (i) $1(A)$ is not in $b$;  (ii) $1(B)$ is not in $b$. 
  In the case (i),  $b'$ is an open branch, and by inductive hypothesis, $v(A) = 0$.  
In the case (ii), $b''$ is an open branch, and by inductive hypothesis, $v(B) = 0$. 
In both cases define $v(A \land  B) = 0$.

\m\mh The valuation $v$  defined in (2.1) and (2.2) clearly satisfies  Definition~\ref{val.sem.letf}.

\mm 

\mh (3) $F = \neg(A \land B)$. 
\m 
 
\mh (3.1)	If $1(\neg(A \land  B))$ is in $b$,  since the  tableau is terminated, \textbf{Rule 3} was applied and the  tableau bifurcated into two branches: $b'$ and $b''$. 
	In $b'$, we have $1(\neg A)$, and in $b''$ $1(\neg B)$. 
	Since the $b$ is an open branch of a terminated tableau, we have two (non-excluding) options: (i) $0(\neg A)$ is not in $b$, (ii) $0(\neg B)$ is not in $b$. 
		If (i), $b'$ is an open branch and, by inductive hypothesis, $v(\neg A) = 1$. 
		Define  $v(\neg(A \land B)) = 1$. 
	If (ii), $b''$ is an open branch and, by inductive hypothesis, $v(\neg B) = 1$. Define $v(\neg(A \land B)) = 1$.
 \m 
 		
	\mh (3.2)	If $0(\neg(A \land  B))$ is in $b$, then, since the  tableau is terminated, \textbf{Rule 4} was applied, therefore $0(\neg A)$ is in $b$ and $0(\neg B)$ is in $b$. But since the  tableau is open, $1(\neg A)$ is not in $b$ and $1(\neg B)$ is not in $b$. Then, by inductive hypothesis, $v(\neg A) = 0$ and $v(\neg B) = 0$. Define  $v(\neg(A \land  B)) = 0$.

\m\mh The valuation $v$  defined  in (3.1) and (3.2) satisfies  Definition~\ref{val.sem.letf}.

	\mm 
	
\mh The cases (4) $F = A \lor B$ and (5) $F = \neg(A \lor B)$ are left to the reader.  	 
	
	\mm
	
\mh (6) $F = \neg\neg A$. 

\m\mh If $1(\neg\neg A)$ is in $b$, by \textbf{Rule 9}, $1(A)$ is also in $b$. 
By inductive hypothesis, 
$v(A)=1$, and we define $v(\neg\neg A) = 1$. By analogous reasoning if $0(\neg\neg A)$  in $b$, 
we define $v(\neg\neg A) = 0$. 

\mh 	The valuation $v$ so defined satisfies Definition~\ref{val.sem.letf}.  
		
	\mmm 
	
\mh We have just shown  that the valuation $v$ defined above satisfies clauses (v1) to (v5) of Definition~\ref{val.sem.letf}. 
It remains to be shown that $v$ also satisfies clauses (v6) and (v7).

\m\mh (7)  $F = \cons A$.

\m\mh (7.1) If $ 1(\cons A) $ is in $b$, then by item (1)  above,  
$v(\cons A) = 1$ and $v(\incon A) = 0$. 
As the tableau is terminated, \textbf{Rule 11} was applied and the  tableau bifurcated into two branches: $b'$ 
and $b''$ such that  
(i) $1(A)$ and $0(\neg A)$ occur in $b'$, and (ii) $0(A)$ and $1(\neg A)$ occur in $b''$. 
If $b'$ is open, then by inductive hypothesis $v(A)=1$ and $v(\neg A)=0$, and  if $b''$ is open, 
then by inductive hypothesis $v(A)=0$ and $v(\neg A)=1$. 
 So, $v$  
satisfies Definition~\ref{val.sem.letf} (clause(v6)).

\mm

\mh (7.2) If $0(\cons A)$ is in $b$, then, as $b$ is open, $1(\cons A)$ is not in $b$, 
and by item (1)  above,  
$v(\cons A) = 0$ and $v(\incon A) = 1$. 
	
	\mm 
	
\mh (8)  $F = \incon A$. 

\m\mh (8.1) If $1(\incon A)$ is in $b$, then, as the  tableau is terminated,  \textbf{Rule 12} was applied and $0(\cons A)$ is in $b$. 
Since the  tableau is open, $1(\cons A)$ is not in $b$. 
By item (1)  above,  $ v(\incon A) = 1$ and $v(\cons A) = 0$. 

\m\mh (8.2) If $0(\incon A)$ is in $b$, then, as the  tableau is terminated, $1(\cons A)$ is in $b$. 
By item (1)  above,  $ v(\incon A) = 1$ and $v(\cons A) = 0$. 
 	
		\mm 
	
\mh (9) $F = \neg\cons A$.  	

\m\mh If $1(\neg\cons A)$ is in $b$, $v(\neg\cons A) = 1$ by definition. 
If $1(\neg\cons A)$ is not in $b$, $v(\neg\cons A) = 0$ by definition. 
 Analogous reasoning applies to  $F = \neg\incon A$.

\mm

\mh Therefore,  $v$ as defined is a legitimate \letf-valuation. 
 	 \end{proof}
	\end{Le}

\m 

\begin{Th} {[Completeness]} \label{completeness.th}
 If $\Gamma \vDash A$, then $\Gamma \vdash A$.
		\begin{proof}
	We prove the contrapositive: if $\Gamma \nvdash A$, then $\Gamma \nvDash A$. 
	Suppose $\Gamma \nvdash A$. Thus there is a terminated \letf-tableau with at least one open branch $b$ such that 
	$0(A)$ is in $b$ and for every formula $F$ from $\Gamma$, $1(F)$ is in $b$. 
	By  Lemma \ref{induced.val},  
	there is a \letf-valuation $v$ induced by $b$ 
	such that 
if	$1(F)$ is in $b$, $v(F) = 1$, and if $0(F)$ is in $b$, $v(F) = 0$. 
	Therefore, there is a \letf-valuation $v$ such that $v(A) = 0$, and for every formula $F$ from $\Gamma$,  $v(F) = 1$. 
	Therefore, $\Gamma \nvDash A$.
		\end{proof}
	
\end{Th}

Clearly, the tableau system  for \letf\ introduced here is
equivalent  to the natural deduction  formulation presented in \cite{letf}.  Indeed,
Theorem  \ref{completeness.th}  shows that  the tableau system for  \letf\  is semantically
characterized by the same valuation semantics (Definition \ref{def.val.semi}) that
characterizes  the  natural deduction rules for  the version of   \letf\ 
introduced in \cite{letf}.

 \subsection{Decidability} \label{subsec.dec}

\begin{Df} \label{gen.subf} {[Generalized  subformula]} 

\begin{enumerate} \setl     
\item  If  \textit{B} is  a subformula of \textit{A} (in the usual sense)  and $A \not = B$, 
then \textit{B} is an {immediate subformula} of \textit{A}. 
\item   If \textit{B} is an immediate subformula of \textit{A}, then \textit{B} is a
{generalized subformula} of \textit{A}.
\item $\neg A$ and $\neg B$ are    generalized subformulas of both  $\neg (A \land B)$    
and $\neg (A \lor B)$.
\item $\neg A$ is a generalized subformula of $\cons A$. 
\item $\cons A$ is a generalized subformula of $\incon A$. 
\item   If \textit{C} is a  generalized subformula of \textit{B} and  \textit{B} is an  generalized
subformula of \textit{A}, then  \textit{C} is a  generalized subformula of \textit{A}.    
\end{enumerate}
 
\end{Df}

\mh As a consequence  of the definition above,  
both $A$ and $\neg A$ are generalized subformulas of $\cons A$  and $\incon A$, since 
$\cons A$ is a generalized  subformula of $\incon A$. 
Besides, in view  of the Definition~\ref{gen.subf},  it is easy to see that if 
\textit{B} is a generalized subformula  of \textit{A},  then $\mathcal{C}(B) <  \mathcal{C}(A)$.
 
\m 

\begin{Th}  {[Decidability]} \

\mh \letf-tableaux provide a decision procedure for \letf. 

\begin{proof} 
Clearly, every term occurring in a  
\letf-tableau of $\Gamma\vdash A$ consists of  signed formulas of $\Gamma \cup \{A\}$ (in the first node) and 
of signed generalized  subformulas of $\Gamma \cup \{A\}$ (in the subsequent nodes),   
and   each~tableau rule yields generalized  subformulas of the formula
to which the rule is applied.  
 Since the  complexity of  formulas occurring in the tableau 
is monotonically decreasing by applications of rules,   
 all tableau branches   are  either closed   
or reach formulas of less complexity 
for which there is no rule to be applied, namely, a  literal (with label 0 or 1), 
$\cons A$ (with label 0), $\neg\cons A$ or $\neg\incon A$ (with label 0 or 1). 
Therefore \letf-tableaux provide a
decision procedure for \letf. 
 
\end{proof} 

\end{Th}

\section{Some examples of \letf-tableaux} \label{examples}

In this section, we give some examples of  tableaux that illustrate properties of \letf. 

\mm 

\begin{Ex} {[Bottom particle]}

\mm

\mh A bottom particle can be defined in \letf\ as $p \land \neg p \land \cons p$, and clearly, $\bot \vdash B$, 
for any \textit{B}.   
%
  
  
  \mm

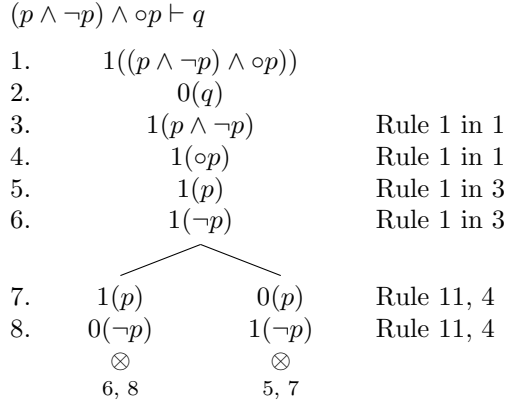
{
	\begin{tableau}
		{
			to prove={ (p \land \neg p) \land \cons p \vdash  q }, 
			for tree={s sep=10mm} 
		}
		[1((p \land \neg p) \land \cons p))
		[0(q)
		[1(p \land \neg p), just={Rule 1 in 1}
		[1(\cons p), just={Rule 1 in 1}
		[1(p), just={Rule 1 in 3}
		[1(\neg p), just={Rule 1 in 3}
		[1(p), just={Rule 11, 4}
		[0(\neg p), just={Rule 11, 4}, close={6, 8}]]
		[0(p)
		[1(\neg p), close={5, 7}]]
		]
		]
		]
		]
		]
		]
	\end{tableau}
}


 \mm 
 
 \end{Ex}
   
   \mm

\mh Every \lfi\ has a bottom particle, since in every \lfi\ a bottom particle can be defined as above.  
This is because the principle of gentle explosion (item \ref{gent.exp} page \pageref{gent.exp}) is an essential feature of \lfis. 

\mm

    \begin{Ex} {[Recovering modus ponens]} \ 
    
    \mm 
   
   \mh As expected, disjunctive syllogism does not hold in \letf, and so modus ponens, since the natural way of defining $A \to B$ 
in \letf\ is  as $\neg A \lor B$. 

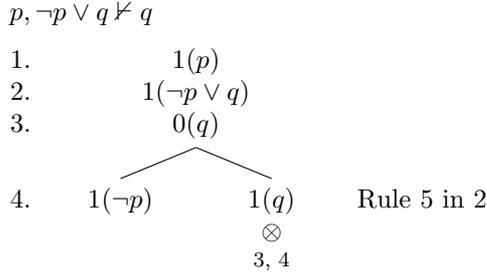
\mm 
     
{
	\begin{tableau}
		{
			to prove={ p ,   \neg p \lor q   \nvdash q },   
			for tree={s sep=10mm} 
		}
	[1(p) 
	[1(\neg p \lor q) 
	[0(q) 
	[1(\neg p), just={\textrm{Rule 5 in 2}}]
	[1(q), close = {3, 4}]
	]
	]
	]
	\end{tableau}
}

 \mm

\mh The open branch gives  a counterexample:  $v(p) = v(\neg p) = 1$ and $v(q)=0$. 
For classical $p$ 
modus ponens is recovered, as we see below \cite[cf.][Fact~32]{letf}.   

\mm 

{
	\begin{tableau}
		{ to prove={ \cons p, p, \neg p  \lor  q  \vdash q}, 
			for tree={s sep=6mm} 
		}
		[1(\cons p)
		[1(p)
		[1(\neg p \lor q)
		[0(q)
				[1(\neg p), just={Rule 5 in 3}
		[1(p), just={Rule 11 in 1}
		[0(\neg p), just={Rule 11 in 1}, close={5, 7}]]
				[0(p), just={Rule 11 in 1}
		[1(\neg p), just={Rule 11 in 1}, close={2, 6}]]]
		[1(q), close={4, 5}]
		]
		]
		]
		]
	\end{tableau}
	
}

\end{Ex}

\mm

\begin{Ex} {[Proofs by cases]} \label{ccc}

\mm
\mh
Excluded middle does not hold in \letf, so the usual form of proof by cases does not obtain. 
But  \letf\ allows other forms of proof  by cases, e.g., from  $\vdash \cons p \lor \incon p$ and 
$\vdash p \lor \neg p \lor \incon p $. 

\mm 

{
	\begin{tableau}
		{to prove={\vdash \cons p \lor \incon  p}, 
			for tree={s sep=10mm}} 
[0(\cons p \lor \incon p )
[0(\cons p), just={Rule 6 in 1}
[0(\incon p), just={Rule 6 in 1}
[1(\cons p), just={Rule 13 in 3}, close={2, 4} ]]]]
	\end{tableau}
	
}

\mm\mh The operators \incon\ and \cons\ work  as if $\incon A$ were a classical negation of  $\cons A$, and vice-versa. 
Indeed, $\cons A, \incon A \vdash  B$ prohibits that $\cons A$ and $\incon A$  hold together, on pain of triviality.  
Note, however, that $\cons A \lor \neg\cons A$ is not valid, as we see below.  \mm 

 {
	\begin{tableau}
		{
			to prove={ \nvdash \cons p \lor \neg\cons p}, 
			for tree={s sep=10mm} 
		}
		[0(\cons p \lor \neg\cons p)
		[0(\cons p), just={Rule 6 in 1}
		[0(\neg\cons p), just={Rule 6 in 1}
		]
		]
		]
	\end{tableau}
}

\mm
\mh 
Formulas $\cons p$ and $\neg\cons p$ are independent 
of each other, as we have seen in the  quasi-matrice of Example  \ref{ex.valuations.letf}.

\end{Ex}

\mm

 The Example \ref{ccc} above indicates  that $\neg\cons p$ and $\incon p$ are not equivalent, as we see below. 

\begin{Ex}

 \mm 

  {
	\begin{tableau}
		{
			to prove={  \neg \cons p \nvdash \incon p}, 
			for tree={s sep=10mm} 
		}
		[1(\neg \cons p)
		[0(\incon p)
		[1(\cons p), just={Rule 13 in 2}
		[1(p), just={Rule 11 in 3}
		[0(\neg p), just={Rule 11 in 3}]]
		[0(p), just={Rule 11 in 3}
		[1(\neg p), just={Rule 11 in 3}]]
		]
		]
		]
	\end{tableau}
}

\mm

\mh The proof of $\incon p \nvdash \neg\cons p$ is left to the reader. 
In \letf, negation of classicality does not entail classicality because the negation \wneg\ is still a weak negation. 
The same applies to $\cons p$ and $\neg\incon p$, which are not equivalent.  
As far as we know, classical negation cannot be defined in \letf\ \cite[cf.][footnote~15]{letf}.

\end{Ex}

\mm
\begin{Ex} {[On `quasi-negations'  in \letf]}

\mh Two unary connectives that are (in some sense) negations can be defined in \letf. 
\begin{enumerate}[label={(\arabic*)}]
\item $\sneg A \vcentcolon  = \cons A \land \neg A$, 
\item  ${\approx} A \vcentcolon  = \incon A \lor \neg A$.\footnote{The connectives \sneg\ and $\approx$ are called, respectively, 
supplementing and complementing negation  \cite[cf.][Def.~33]{letf}. }  
\eenr

\mh We have that explosion holds for \sneg\ and excluded middle for $\approx$, i.e.   
\enr 
\item $A, \sneg A \vdash B$, 
\item  $\vdash A \lor {\approx} A $,   
\eenr
\mh as well as the following dual inferences, 
\enr 
\item $\sneg A \vdash \neg A$, 
\item  $\neg A \vdash {\approx} A $.    
\eenr
\mh 
On the other hand,  neither double negation nor de Morgan 
  hold for   \sneg\ and $\approx$ in \letf. 
  For this reason, we think these connectives  should rather be called 
  `quasi-negations'. We prove below 
$p \nvdash {\approx}{\approx}  p$ and $\sneg (p \lor q) \nvdash \sneg p \land \sneg q$, 
and leave to the reader the other invalid inferences  with \sneg\ and ${\approx}$.

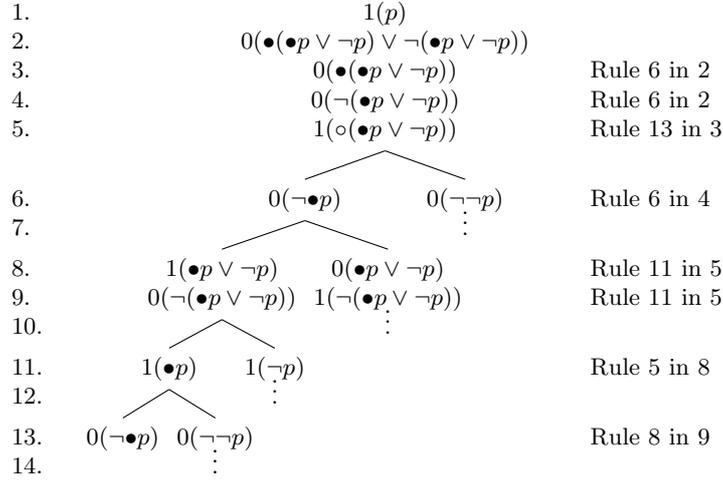
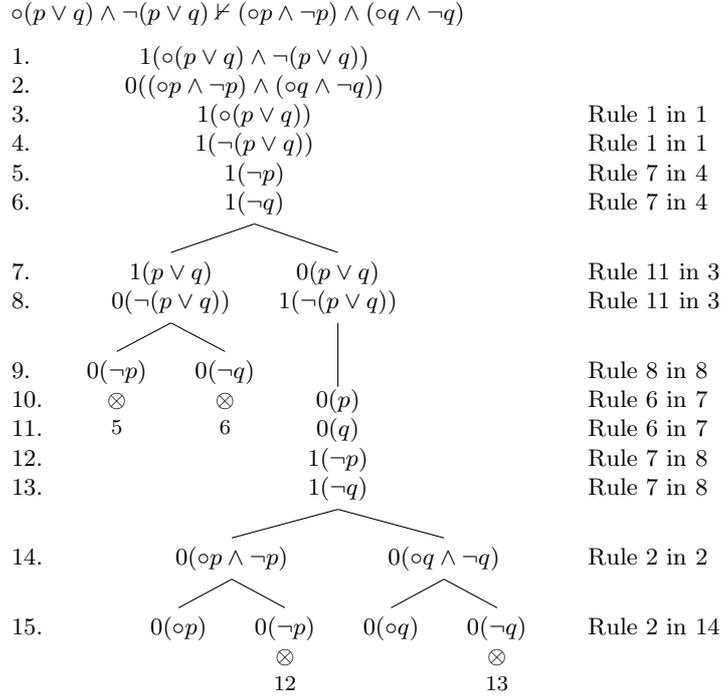
\mm 

{\footnotesize 
	\begin{tableau}
		{to prove={p \nvdash \incon (\incon p \lor \neg p) \lor \neg (\incon p \lor \neg p)}, 
			for tree={s sep=0mm} 
		}
		[1(p) 
		[0(\incon (\incon p \lor \neg p) \lor \neg (\incon p \lor \neg p))
		[0(\incon (\incon p \lor \neg p)), just={Rule 6 in 2}   
		[0(\neg (\incon p \lor \neg p)), just={Rule 6 in 2}     
[1(\cons (\incon p \lor \neg p)), just={Rule 13 in 3}		
		[0(\neg\incon p) 
		[1(\incon p \lor \neg p), just={Rule 11 in 5} 
		[0(\neg(\incon p \lor \neg p)), just={Rule 11 in 5} 
		[1(\incon p) [0(\neg\incon p)] 
		[0(\neg\neg p), just={Rule 8 in 9} 
		[ \vdots  
		]]] 
		[1(\neg p), just={Rule 5 in 8} [\vdots]]  ]		
		] 
		[0(\incon p \lor \neg p) [1(\neg(\incon p \lor \neg p)) [\vdots]] ]     ]
		[0(\neg\neg p), just={Rule 6 in 4} [ \vdots   
		]]  		
		]
		] 		
		]]
		] 
		]]
		]
		]
		]
	\end{tableau}
	
}

\mm

{\footnotesize

	\begin{tableau} 
		{to prove={\cons (p \lor q) \land \neg(p \lor q) \nvdash (\cons p \land\neg p) \land (\cons q \land \neg q)}, 
			for tree={s sep=4mm} 
		} 
		[1(\cons (p \lor q) \land \neg(p \lor q))
		[0((\cons p \land\neg p) \land (\cons q \land \neg q))
		[1(\cons(p \lor q)), just={Rule 1 in 1}
		[1(\neg(p \lor q)), just={Rule 1 in 1} 
		[1(\neg p), just={Rule 7 in 4}
		[1(\neg q), just={Rule 7 in 4} 
		[1(p \lor q), just={Rule 11 in 3}
		[0(\neg(p \lor q)), just={Rule 11 in 3}  
[0(\neg p), close={5}] [0(\neg q), just={Rule 8 in 8 }, close={6}]	
		]]
		[0(p \lor q)
		[1(\neg(p \lor q)) 
[0(p), move by=1 , just={Rule 6 in 7 }
		[0(q), just={Rule 6 in 7 }		
		[1(\neg p), just={Rule 7 in 8 } 
		[1(\neg q), just={Rule 7 in 8 }
				[0(\cons p \land \neg p)  [0(\cons p)] [0(\neg p), just={Rule 2 in 14} , close={12}]
		] 
		[0(\cons q \land\neg q), just={Rule 2 in 2}   [0(\cons q)] [0(\neg q), close={13}] ]
		]]     ]] ] 		] 
		]]] ] 		]]]
		]		]]		]
	\end{tableau}  
	}

\end{Ex}

     \mmm

 A central point of the intended interpretation of \lets\  is that there may be scenarios in which 
there is evidence for exactly one among 
$A$ and $\neg A$, but   such evidence is still not conclusive. In this case, $A \lor \neg A$ holds, as well 
$\neg(A \land\neg A)$, but $\cons A$  does not hold. 
Indeed, in a  scenario such that $A$ holds but $\neg A$ does not hold, the evidence   for $A$ may   be non-conclusive, 
and so $A$ cannot be taken as true and is not subjected to classical logic. 
This idea is expressed in the semantics by means of clause (v6), which states only a necessary condition for $\cons A$. 
The same idea could be expressed equivalently, with \incon,  by the contrapositive of clause (v6) -- indeed, 
when both  $A$ and  $\neg A$ hold, as well as when neither holds, $\incon A$ holds.  
The behavior of \cons\ and \incon\ 
is illustrated by   Example \ref{ex.x} below.

\begin{Ex} {[On the behavior of \cons\ and \incon]} 
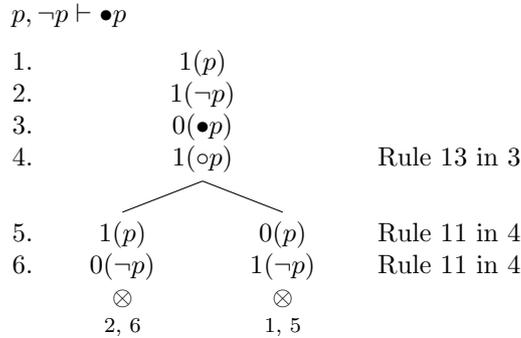
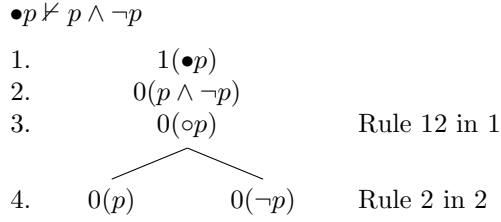
\label{ex.x}

\m

{
	\begin{tableau}
		{
			to prove={   p, \neg p \vdash \incon p}, 
			for tree={s sep=10mm} 
		}
		[1(p)
		[1(\neg p)
		[0(\incon p)
		[1(\cons p), just={Rule 13 in 3}
		[1(p), just={Rule 11 in 4}
		[0(\neg p), just={Rule 11 in 4}, close={2, 6}]]
		[0(p), just={Rule 11 in 4}
		[1(\neg p), just={Rule 11 in 4}, close={1, 5}]]
		]
		]
		]
		]
	\end{tableau}
}

\mm

   {
	\begin{tableau}
		{
			to prove={  \incon p \nvdash p \land \neg p}, 
			for tree={s sep=10mm} 
		}
		[1(\incon p)
		[0(p \land \neg p)
		[0(\cons p), just={Rule 12 in 1}
		[0(p), just={Rule 2 in 2}]
		[0(\neg p), just={Rule 2 in 2}]
		]
		]
		]
	\end{tableau}
}

\mm

\mh Proofs of $\cons p \vdash p \lor\neg p$ and $p \lor\neg p \nvdash \cons p$ are left to the reader. 

\end{Ex}

\begin{Ex} {[Propagation of classicality]}

\mh The operator \cons\ does not propagate 
over more complex formulas, that is, $\cons A , \cons B $ does not imply $\cons \neg A$, $\cons (A \land B)$, or 
$\cons (A \lor B)$.  We illustrate this fact showing that $\cons p , \cons q \nvdash \cons (p \land q)$. 

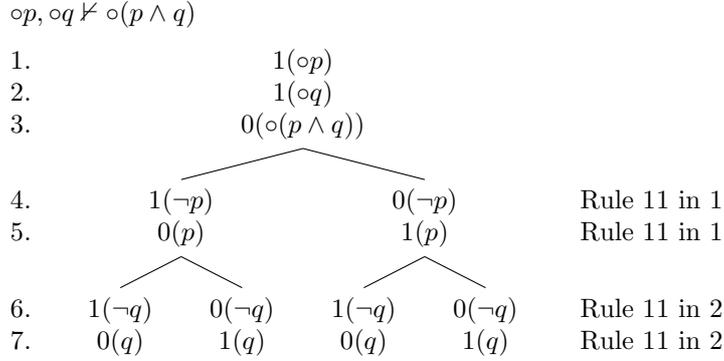
\mm

{
	\begin{tableau}
		{to prove={\cons p , \cons q \nvdash \cons (p \land q)}, 
			for tree={s sep=5mm} 
		}
		[1(\cons p)
		[1(\cons q)
		[0(\cons (p \land  q))
		[1(\neg p), just={Rule 11 in 1}
		[0(p), just={Rule 11 in 1} 
        [1(\neg q), just={Rule 11 in 2} [0(q), just={Rule 11 in 2}]
        ] 
        [0(\neg q) [1(q)]
        ] 		
		]]
		[0(\neg p)
		[1(p)
		[1(\neg q)  [0(q)]
        ] 
        [0(\neg q) [1(q)]
        ] 		
		]] 		
		]]
		] 
		]]
		]
		]
		]
	\end{tableau}
	
}

\mm
\mh Notice, however, that the set $\{  \cons p, \cons q, p \land q, \neg(p \land q) \}$ is not satisfiable.   
Indeed, although in \letf\ \cons\ does not propagate over \wneg,  $\land$, and $\lor$, 
 if  $\cons p_1, ... , \cons p_n$ hold, then    all formulas formed with $\{p_1, ..., p_n  \}$ 
over $\{ \neg, \land, \lor \}$ behave classically, as has been shown in \cite[Fact~31]{letf}.  
\end{Ex}

\section{Final remarks}

Analytic  tableaux constitute a decision procedure for \letf\ that, we think, is at least more  elegant  than quasi-matrices. However,  besides issues  of  elegance, 
it is reasonable to conjecture that \letf-tableaux are  in fact 
more efficient than the quasi-matrices.

D'Agostino showed in  \cite{dagostino} that classical tableaux are globally less efficient than truth-tables as a decision 
procedure for classical logic. 
However, this result certainly does not apply to any logic, and it would perhaps have to be revised in the case of some non-classical logics. 
Non-deterministic matrices, like the ones yielded by the semantics of \letf\ presented in \cite{letf}, are considerably greater in size than 
classical  truth-tables, but at first sight the same does not occur with respect to \letf-tableaux when compared to classical tableaux. 
In the case of \letf, the examples of Section \ref{examples}  suggest  that, at least locally, for some classes of formulas, tableaux are 
indeed better than quasi-matrices.  For this reason, whether or not D'Agostino's results apply  
to \letf \  is a question that deserves to be further investigated. 

\letf\ is an extension of \fde, or 
Belnap-Dunn logic,   
a  widely studied four-valued  paraconsistent logic  
with applications in computer science  and in the algebra of processes. 
\letf\  extends \fde\   by adding  a classicality operator \cons\ and  its dual,
 a non-classicality operator  \incon.  Propositions $\cons A$ and $\incon A$ can be read,  
 respectively, as saying that the evidence for $A$ is conclusive and non-conclusive, as well as that 
 the information about $A$ is reliable and non-reliable. 
Considering its  potential
applications  in automated reasoning and artificial intelligence, \letf\ 
deserves in-depth investigations. 
This paper is one more step
in this direction.

\bibliographystyle{plainnat}       
\bibliography{refs.answer}

\begin{thebibliography}{19}
\providecommand{\natexlab}[1]{#1}
\providecommand{\url}[1]{\texttt{#1}}
\expandafter\ifx\csname urlstyle\endcsname\relax
  \providecommand{\doi}[1]{doi: #1}\else
  \providecommand{\doi}{doi: \begingroup \urlstyle{rm}\Url}\fi

\bibitem[Antunes et~al.(2020)Antunes, Carnielli, Kapsner, and
  Rodrigues]{kpk.let}
H.~Antunes, W.~Carnielli, A.~Kapsner, and A.~Rodrigues.
\newblock Kripke-style models for logics of evidence and truth.
\newblock \emph{Axioms}, 9, 2020.

\bibitem[Belnap(1977)]{belnap1977.how}
N.~D. Belnap.
\newblock How a computer should think.
\newblock In G.~Ryle, editor, \emph{Contemporary {A}spects of {P}hilosophy}.
  Oriel Press \cite[reprinted in][]{essays.belnap.dunn}, 1977.

\bibitem[Bergstra and Ponse(2000)]{bergstra}
A.~Bergstra and A.~Ponse.
\newblock Process algebra with four-valued logic.
\newblock \emph{Journal of Applied Non-Classical Logics}, 10:\penalty0 27--53,
  2000.

\bibitem[Bueno-Soler et~al.()Bueno-Soler, Carnielli, Coniglio, and
  Rodrigues]{future}
J.~Bueno-Soler, W.~Carnielli, M.~Coniglio, and A.~Rodrigues.
\newblock Applications of logics of evidence and truth to automated reasoning.
\newblock \emph{Manuscript, to appear}.

\bibitem[Carnielli and Marcos(2001)]{tableau.lfi}
W.~Carnielli and J.~Marcos.
\newblock Tableau systems for logics of formal inconsistency.
\newblock In \emph{Proceedings of the International Conference on Artificial
  Intelligence (IC-AI 2001)}, pages 848--852, 2001.

\bibitem[Carnielli and Rodrigues(2017)]{letj}
W.~Carnielli and A.~Rodrigues.
\newblock An epistemic approach to paraconsistency: a logic of evidence and
  truth.
\newblock \emph{Synthese}, 196:\penalty0 3789--3813, 2017.

\bibitem[Carnielli et~al.(2007)Carnielli, Coniglio, and Marcos]{carcoma2007}
W.~Carnielli, M.~Coniglio, and J.~Marcos.
\newblock Logics of formal inconsistency.
\newblock In Gabbay~\& Guenthner, editor, \emph{Handbook of Philosophical
  Logic}, volume~14. Springer, 2007.

\bibitem[Carnielli et~al.(2019)Carnielli, Coniglio, and Rodrigues]{recovery}
W.~Carnielli, M.~Coniglio, and A.~Rodrigues.
\newblock Recovery operators, paraconsistency and duality.
\newblock \emph{Logic Journal of the IGPL}, 2019.
\newblock \doi{10.1093/jigpal/jzy054}.

\bibitem[D'{A}gostino(1992)]{dagostino}
M.~D'{A}gostino.
\newblock Are tableaux an improvement on truth-tables?
\newblock \emph{Journal of Logic, Language, and Information}, 1:\penalty0
  235--252, 1992.

\bibitem[Dunn(1976)]{dunn76}
J.~M. Dunn.
\newblock Semantics for first-degree entailments and `coupled trees'.
\newblock \emph{Philosophical Studies}, 29\penalty0 (3):\penalty0 149--168
  \cite[reprinted in][]{essays.belnap.dunn}, 1976.

\bibitem[Ginsberg(1988)]{ginsberg}
M.~L. Ginsberg.
\newblock Multi-valued logics: A uniform approach to reasoning in artificial
  intelligence.
\newblock \emph{Computer Intelligence}, 4:\penalty0 256--316, 1988.

\bibitem[Loparic(1986)]{loparic1986}
A.~Loparic.
\newblock A semantical study of some propositional calculi.
\newblock \emph{The Journal of Non-Classical Logic}, 3\penalty0 (1):\penalty0
  73--95, 1986.

\bibitem[Loparic(2010)]{loparic2010}
A.~Loparic.
\newblock Valuation semantics for intuitionistic propositional calculus and
  some of its subcalculi.
\newblock \emph{Principia}, 14\penalty0 (1):\penalty0 125--133, 2010.

\bibitem[Marcos(2005)]{marcos2005}
J.~Marcos.
\newblock Nearly every normal modal logic is paranormal.
\newblock \emph{Logique et Analyse}, 48:\penalty0 279--300, 2005.

\bibitem[Omori and Wansing(2019)]{essays.belnap.dunn}
H.~Omori and H.~Wansing, editors.
\newblock \emph{New Essays on Belnap-Dunn Logic -- Synthese Library Studies in
  Epistemology, Logic, Methodology, and Philosophy of Science, Volume 418}.
\newblock Springer, 2019.

\bibitem[Patel-Schneider(1989)]{patel}
P.~F. Patel-Schneider.
\newblock A four-valued semantics for terminological logics.
\newblock \emph{Artificial Intelligence}, 38:\penalty0 319--351, 1989.

\bibitem[Rodrigues and Carnielli(to appear)]{answer.barrio}
A.~Rodrigues and W.~Carnielli.
\newblock {On Barrio, Lo Guercio, and Szmuc on logics of evidence and truth}.
\newblock \emph{Logic and Logical Philosophy}, to appear.

\bibitem[Rodrigues et~al.(2020)Rodrigues, Bueno-Soler, and Carnielli]{letf}
A.~Rodrigues, J.~Bueno-Soler, and W.~Carnielli.
\newblock Measuring evidence: a probabilistic approach to an extension of
  {B}elnap-{D}unn logic.
\newblock \emph{Synthese}, 2020.
\newblock \doi{10.1007/s11229-020-02571}.

\bibitem[Straccia(1997)]{straccia}
U.~Straccia.
\newblock A sequent calculus for reasoning in four-valued description logics.
\newblock In D.~Galmiche, editor, \emph{Automated Reasoning with Analytic
  Tableaux and Related Methods}. Berlin, Springer, 1997.

\end{thebibliography}

\end{document}